# Analytic continuation of representations and estimates of automorphic forms

By Joseph Bernstein and Andre Reznikov

### 0. Introduction

0.1. *Analytic vectors and their analytic continuation.* Let $G$ be a Lie group and $(\pi, G, V)$ a continuous representation of $G$ in a topological vector space $V$. A vector $v \in V$ is called *analytic* if the function $\xi_v : g \mapsto \pi(g)v$ is a real analytic function on $G$ with values in $V$. This means that there exists a neighborhood $U$ of $G$ in its complexification $G_{\mathbb{C}}$ such that $\xi_v$ extends to a holomorphic function on $U$. In other words, for each element $g \in U$ we can unambiguously define the vector $\pi(g)v$ as $\xi_v(g)$, i.e., we can extend the action of $G$ to a somewhat larger set. In this paper we will show that the possibility of such an extension sometimes allows one to prove some highly nontrivial estimates.

Unless otherwise stated, $G = \mathrm{SL}(2, \mathbb{R})$, so $G_{\mathbb{C}} = \mathrm{SL}(2, \mathbb{C})$. We consider a typical representation of $G$, i.e., a representation of the principal series. Namely, fix $\lambda \in \mathbb{C}$ and consider the space $\mathfrak{D}_\lambda$ of smooth homogeneous functions of degree $\lambda - 1$ on $\mathbb{R}^2 \setminus 0$, i.e., $\mathfrak{D}_\lambda = \{\phi \in C^\infty(\mathbb{R}^2 \setminus 0) : \phi(ax, ay) = |a|^{\lambda-1}\phi(x, y)\}$; we denote by $(\pi_\lambda, G, \mathfrak{D}_\lambda)$ the natural representation of $G$ in the space $\mathfrak{D}_\lambda$.

Restriction to $S^1$ gives an isomorphism $\mathfrak{D}_\lambda \simeq C^\infty_{\text{even}}(S^1)$, and for basis vectors of $\mathfrak{D}_\lambda$ one can take the vectors $e_k = \exp(2ik\theta)$. If $\lambda = it$, then $(\pi_\lambda, \mathfrak{D}_\lambda)$ is a unitary representation of $G$ with the invariant norm $||\phi||^2 = \frac{1}{2\pi}\int_{S^1} |\phi|^2 d\theta$.

Consider the vector $v = e_0 \in \mathfrak{D}_\lambda$. We claim that $v$ is an analytic vector and we want to exhibit a large set of elements $g \in G_{\mathbb{C}}$ for which the expression $\pi(g)v$ makes sense. The vector $v$ is represented by the function $(x^2 + y^2)^{\frac{\lambda-1}{2}} \in \mathfrak{D}_\lambda$. For any $a > 0$ consider the diagonal matrix $g_a = \mathrm{diag}(a^{-1}, a)$. Then

$$\xi_v(g_a) = \pi_\lambda(g_a)v = (a^2 x^2 + a^{-2} y^2)^{\frac{\lambda-1}{2}}.$$

This last expression makes sense as a vector in $\mathfrak{D}_\lambda$ for any *complex* $a$ such that $|\arg(a)| < \frac{\pi}{4}$ (since in this case $\mathrm{Re}\,(a^2 x^2 + a^{-2} y^2) > 0$). Hence, we see that the function $\xi_v$ extends analytically to the subset $I = \{g_a : |\arg(a)| < \frac{\pi}{4}\} \subset \mathrm{SL}(2, \mathbb{C})$.

The same argument shows that the function $\xi_v$ extends analytically to the domain $U = \mathrm{SL}(2, \mathbb{R}) \cdot I \cdot K_{\mathbb{C}} \subset \mathrm{SL}(2, \mathbb{C})$ (open in the usual topology), where $K = \mathrm{SO}(2,\mathbb{R})$ and $K_{\mathbb{C}} = \mathrm{SO}(2,\mathbb{C}) \simeq \mathbb{C}^*$; thus, for any $g \in U$ we unambiguously define the vector $\pi(g)v$.

As $g$ approaches the boundary of $U$, the vector $\pi(g)v \in \mathfrak{D}_\lambda$ has very specific asymptotic behavior that we will use in order to obtain information about this vector.

0.2. *Triple products.* Let us describe an application of the principle of analytic continuation to a problem in the theory of automorphic functions. Namely, we will show how to apply the principle in order to settle a conjecture of Peter Sarnak on triple products. As a corollary of our result we will get a new bound on Fourier coefficients of cusp forms.

Recall the setting. Let $\mathfrak{H}$ be the upper half-plane with the hyperbolic metric of constant curvature $-1$. We consider the natural action of the group $G = \mathrm{SL}(2,\mathbb{R})$ on $\mathfrak{H}$ and identify $\mathfrak{H}$ with $G/K$ by means of this action.

Fix a lattice $\Gamma \subset G$ and consider the Riemann surface $Y = \Gamma \backslash \mathfrak{H}$. In this paper we will discuss both cocompact and noncocompact lattices of finite covolume. For simplicity of exposition, in most of the paper we will only discuss the cocompact case. Then in Section 4 we will describe how to overcome the extra difficulties in case of noncocompact lattices.

The Laplace-Beltrami operator $\Delta$ acts on the space of functions on Y. When $Y$ is compact it has discrete spectrum; we denote by $\mu_0 < \mu_1 \leq \ldots$ its eigenvalues on $Y$ and by $\phi_i$ the corresponding eigenfunctions. (We assume that $\phi_i$ are $L^2$ normalized: $||\phi_i|| = 1$.) These functions $\phi_i$ are usually called *automorphic functions* or *Maass forms* (see [B]).

To state the problem about triple products, fix one automorphic function, $\phi$, and consider the function $\phi^2$ on $Y$. Since $\phi^2$ is not an eigenfunction, it is *not* an automorphic function. Since $\phi^2 \in L^2(Y)$, we may consider its spectral decomposition in the basis $\{\phi_i\}$:

$$\phi^2 = \sum c_i \phi_i.$$

Here the coefficients are given by the triple product integrals: $c_i = \langle \phi^2, \phi_i \rangle = \int_X \phi \cdot \phi \cdot \overline{\phi}_i dx$. Later we will explain why these triple products are of interest and how they are related to the theory of Rankin-Selberg $L$-functions (see also [S], which was our starting point).

CLAIM. *The coefficients $c_i$ decay exponentially as* $\exp(-\frac{\pi}{2}\sqrt{\mu_i})$.

More precisely, let us introduce new parameters $\lambda_i$ such that $\mu_i = \frac{1-\lambda_i^2}{4}$ (the meaning of this parametrization will become clear in subsection 0.3). Introduce new (normalized) coefficients $b_i = |c_i|^2 \exp(\frac{\pi}{2}|\lambda_i|)$. The main result of the paper is the proof of the following theorem which settles a conjecture of P. Sarnak (see [S]):

THEOREM. *There exists a constant $C > 0$ such that*

$$\sum_{|\lambda_i| \leq T} b_i \leq C \cdot (\ln T)^3 \quad \text{as } T \to \infty.$$

COROLLARY. *There exists a constant $C > 0$ such that*

$$\left| \int_X \phi^2 \cdot \overline{\phi}_i dx \right| \leq C(\ln \mu_i)^{\frac{3}{2}} \cdot \exp(-\frac{\pi}{2}\sqrt{\mu_i}).$$

*Remarks.* 1. The bound in the theorem is essentially sharp. Namely, our method gives the following lower bound on the average: $\sum_{i=1}^{\infty} b_i e^{-\varepsilon|\lambda_i|} \geq c|\ln \varepsilon|$.

For a single triple product we cannot do better than the bound in the corollary.

For congruence subgroups we can speculate about the true "size" of these triple products. It is known (see 0.6) that in certain cases the $c_i$ are equal (up to an explicit factor) to the value of the triple Garrett $L$-function at $\frac{1}{2}$. For these $L$-functions, the Lindelöf conjecture predicts $b_i \ll |\lambda_i|^{-2+\varepsilon}$. This is consistent with our bound together with the Weyl law: the number of eigenfunctions with $|\lambda_i| \leq T$ is proportional to $T^2$.

2. We will prove similar results for nonuniform lattices (see §4).

3. This type of question has been considered before. The first result on exponential decay of the coefficients $c_i$ for a holomorphic cusp form $\phi$ was proven by A. Good ([G]) for the general (i.e., nonarithmetic) nonuniform lattices $\Gamma$ thanks to a special feature of holomorphic Poincaré series. Recently, M. Jutila ([J]) extended these results to the nonholomorphic case (Maass forms), but only for the group $SL(2, \mathbb{Z})$, using Kuznetsov's formula and nontrivial arithmetic information (Weil's bounds on Kloosterman's sums and deep results of Iwaniec). In particular, all these methods work only for nonuniform lattices.

In [S], P. Sarnak introduced a new method to estimate the triple products based on analytic continuation of certain matrix coefficients of the function $\phi$; this method works for uniform lattices as well. Being partly based on the theory of spherical harmonics, it led to a weaker bound (by a power of $T$).

Our method, in addition to the analytic continuation, uses more sophisticated representation theory, in particular, an idea of $G$-invariant norms on representations and gives the optimal result (possibly, up to a power of logarithm).

4. Our method gives a more general result than Theorem 0.2. We can obtain similar logarithmic bounds for any polynomial expression in any finite number of automorphic functions $\phi_k$ instead of $\phi^2$, as above.

5. One can ask the same question about growth of triple products for polynomial expressions in automorphic functions of nonzero weight. In this case the decay is also exponential with the same exponent as in Claim 0.2, but the bound in the analogue of Theorem 0.2 is *a power* of $T$ and *not* logarithmic as above.

The main interest in triple products and their bounds stems from their relation to the theory of automorphic $L$-functions. We will discuss this relation in 0.6. We also show in 0.7 that Theorem 0.2 implies a new bound on the Fourier coefficients of automorphic functions in the case of nonuniform lattices.

0.3. *Automorphic representations.* To explain our method, we first recall the relation of automorphic functions to automorphic representations of $G$.

For a given lattice $\Gamma$ in $G$ we denote by $X$ the quotient space $X = \Gamma \backslash G$. The group $G$ acts on $X$, hence, on the space of functions on $X$. We can identify $\mathfrak{H}$ with $G/K$. Then the Riemann surface $Y = \Gamma \backslash \mathfrak{H}$ is identified with $X/K$. This induces an isometric embedding $L^2(Y) \subset L^2(X)$, the image consisting of all $K$-invariant functions.

For any eigenfunction $\phi$ of the Laplace operator $\Delta$ on $Y$ we consider the closed $G$-invariant subspace $L_\phi \subset L^2(X)$ generated by $\phi$ under the $G$-action. It is known that $(\pi, L) = (\pi_\phi, L_\phi)$ is an irreducible unitary representation of $G$ (see [G6]).

Conversely, fix an irreducible unitary representation $(\pi, L)$ of the group $G$ and a $K$-fixed unit vector $v_0 \in L$. Then any $G$-morphism $\nu : L \to L^2(X)$ defines an eigenfunction $\phi = \nu(v_0)$ of $\Delta$ on $Y$; if $\nu$ is an isometric embedding, then $||\phi|| = 1$. Thus, the eigenfunctions $\phi$ correspond to the tuples $(\pi, L, v_0, \nu)$.

Usually it is more convenient to work with smooth vectors. Let $V = L^\infty$ be the subspace of smooth vectors in $L$. Then $\nu$ gives a morphism $\nu : V \to (L^2(X))^\infty \subset C^\infty(X)$. If $X$ is compact, then $\mathrm{Mor}_G(L, L^2(X)) \simeq \mathrm{Mor}_G(V, C^\infty(X))$. Thus, the eigenfunctions correspond to the tuples $(\pi, V, v_0, \nu : V \to C^\infty(X))$.

All irreducible unitary representations of $G$ with $K$-fixed vector are classified: these are representations of the principal and complementary series and the trivial representation. For simplicity, consider representations of the principal series only. In this case the representation $(\pi, V)$ in the space of smooth vectors is isomorphic to the representation $(\pi_\lambda, \mathfrak{D}_\lambda)$ for some $\lambda = it$ (see 0.1). The eigenvalue of the corresponding automorphic function equals $\mu = \frac{1-\lambda^2}{4}$.

0.4. *The method.* We describe here the idea behind the proof of Theorem 0.2.

Let $L_i \subset L^2(X)$ be the space corresponding to the automorphic function $\phi_i$ as above (see 0.3). Let $\mathrm{pr}_i : L^2(X) \to L_i$ be the orthogonal projection.

Since the function $\phi^2$ is $K$-invariant and there is at most one $K$-fixed vector in each irreducible representation of $SL(2, \mathbb{R})$, we have $\mathrm{pr}_i(\phi^2) = c_i \phi_i$.

Since the $G$-action commutes with the multiplication of functions on $X$,
$$\mathrm{pr}_i((\pi(g)\phi)^2) = \mathrm{pr}_i(\pi(g)(\phi^2)) = c_i \pi_i(g)\phi_i \ .$$

By the principle of analytic continuation, the same identity holds for the complex points $g \in U$ (see 0.1). Since all the spaces $L_i$ are orthogonal, we get the following basic relation for the complex points $g$:

$$(0.4.1) \qquad ||(\pi(g)\phi)^2||^2 = \sum_i |c_i|^2 ||\pi_i(g)\phi_i||^2 \ .$$

Here $||\cdot|| = ||\cdot||_{L^2}$ denotes the $L^2$-norm in $L^2(X)$.

It is important that in (0.4.1) we deal with complex points $g$ and for such $g$ the operators $\pi(g)$ are *nonunitary*. As a result, relation (0.4.1) gives nontrivial information.

Now, consider the behavior of the function $(\pi(g)\phi)^2$ near the boundary of $U$. Take $\varepsilon > 0$ and an element $g_\varepsilon \in U$ which is approximately at the distance $\varepsilon$ from the boundary of $U$. For example, set $g_\varepsilon = \mathrm{diag}(a_\varepsilon^{-1}, a_\varepsilon)$ for $a_\varepsilon = \exp((\frac{\pi}{4} - \varepsilon)i)$.

With shorthand notation, $v_\varepsilon = \pi(g_\varepsilon)e_0$ and $\phi_\varepsilon = \nu(v_\varepsilon)$, formula (0.4.1) becomes

$$(0.4.2) \qquad ||\phi_\varepsilon^2||^2 = \sum |c_i|^2 ||\phi_{i,\varepsilon}||^2 \ .$$

Our goal is to give an *upper* bound on the left-hand side of (0.4.2) and a *lower* bound of each of the $||\phi_{i,\varepsilon}||^2$ as $i \to \infty$ and $\varepsilon \to 0$. The latter problem is simpler since it is invariantly defined in terms of representation theory; thus it can be computed in any model of the representation $\pi_i$ (e.g., in $\mathfrak{D}_{\lambda_i}$). A direct computation gives

$$||\phi_{i,\varepsilon}||^2 \geq c \cdot \exp((\frac{\pi}{2} - \varepsilon)|\lambda_i|) \text{ for some } c > 0.$$

On the other hand, we will prove the bound $||\phi_\varepsilon^2|| \ll |\ln \varepsilon|^3$. These two bounds easily imply Theorem 0.2 (see 2.3).

The last bound follows from the bound $|\phi_\varepsilon(x)| \leq C|\ln \varepsilon|$ which holds *pointwise* on $X$ and which we consider to be our main achievement in this paper. Its proof is based on the use of invariant norms which we now explain.

0.5. *Invariant norms.* The most difficult part of the proof is that of the pointwise bound $|\phi_\varepsilon| \leq C|\ln \varepsilon|$. Note that the $L^2$-norm of $\phi_\varepsilon$ is of order $|\ln \varepsilon|^{\frac{1}{2}}$; hence, the pointwise bound only differs from it by a power of logarithm.

In order to obtain such a bound, we use invariant (non-Hermitian!) norms on the representation $\pi$. Namely, as we have explained, any automorphic function gives rise to an embedding $\nu : \mathfrak{D}_\lambda \to C^\infty(X)$. We consider the supremum norm $N_{\sup}$ on $\mathfrak{D}_\lambda$ induced by $\nu$:

$$N_{\sup}(v) = \sup_{x \in X} |\nu(v)(x)|.$$

For a discussion of $L_p$-norms on $X$ see Appendix A.

From the Sobolev restriction theorem on $X$ (for more details see Appendix B), it follows that $N_{\sup}$ is bounded by some Sobolev norm $S = S_k$ on the space $\mathfrak{D}_\lambda$. Hence, the main properties of the norm $N_{\sup}$ are: it is $G$-invariant and $N_{\sup} \leq S$.

We will show that there exists a maximal norm $S^G$ on the space $\mathfrak{D}_\lambda$ satisfying these two conditions. This norm is defined in terms of the representation $\pi_\lambda$ only and it is independent of the automorphic form picture.

We then use the model $\mathfrak{D}_\lambda$ of $\pi_\lambda$ in order to prove the bound $S^G(v_\varepsilon) \leq C|\ln \varepsilon|$. The proof uses the standard method of dyadic decomposition from harmonic analysis; it is based on the observation that, in $\mathfrak{D}_\lambda$, the vector $v_\varepsilon$ is represented by a function which is roughly homogeneous.

As a result we get a pointwise bound

(0.5) $\qquad \sup |\phi_\varepsilon| = N_{\sup}(v_\varepsilon) \leq S^G(v_\varepsilon) \leq C|\ln \varepsilon|$ as $\varepsilon \to 0$.

*Remark.* A new feature of our method, which seems to be absent in the classical approaches to automorphic forms, is the essential use of representation theory.

First of all, in order to study the automorphic function $\phi$ that lives on the space $Y$, we pass to a bigger space, $X$, and work directly with the representation $(\pi, G, V) \subset C^\infty(X)$ which corresponds to $\phi$.

In some classical approaches, the space $V$ is actually also present, albeit very implicitly. And when present, it appears only as a collection of vectors $\pi(h)\phi$ created from the automorphic function $\phi$ by operators $\pi(h)$ corresponding to various functions (or distributions) $h$ on $G$. Though, in principle, one can show that such functions exhaust $V$, in most cases it is very difficult to work with such an implicit description.

In this paper we directly use the space $V$ in order to prove Theorem 0.2. For example, the central technical result is the pointwise bound of the function $u = \phi_\varepsilon \in V$. This bound is proven in Section 5 by means of dyadic decomposition. The idea of the method is to break the function $u$ into the sum of "pieces" $u_i \in V$ which we can move to a better position (for more detail see §§5.2).

We describe these $u_i$ using the explicit model $\mathfrak{D}_\lambda$ of $V$. We do not know how to realize the $u_i$'s in the form $\pi(h)\phi$. So we do not see how to prove this crucial estimate without using the space $V$ as a whole.

0.6. *Relation to L-functions.* The main interest in triple products and their bounds stems from their relation to the theory of automorphic $L$-functions. A particular case of these triple products is the scalar product of $\phi^2$ with the Eisenstein series $E(s)$. This is the original example of Rankin and Selberg of the $L$-function associated to two cusp forms (see [B]). Namely, $L(\phi \otimes \phi, s) = g(s)\langle \phi^2, E(s)\rangle$, where $g(s)$ is an explicit factor.

M. Harris and S. Kudla ([HK]) discovered that such triple products are related to the special value at $s = \frac{1}{2}$ of $L(\phi \otimes \phi \otimes \phi_i, s)$. This gives further reason for the study of such triple products, at least when $\phi$ and $\phi_i$ are holomorphic cusp forms for a congruence subgroup of a division algebra.

0.7. *Bounds on Fourier coefficients of cusp forms.* As we mentioned above, our result implies certain bounds for the Rankin-Selberg $L$-functions on the critical line. This, in turn, has implication for the classical problem of obtaining bounds of the Fourier coefficients of cusp forms.

Recall the setting (see [Se], [G], [S]). Let $\Gamma$ be a nonuniform lattice in $\mathrm{SL}(2,\mathbb{R})$, which can be nonarithmetic (the standard example of a nonuniform lattice is $\Gamma = \mathrm{SL}(2,\mathbb{Z})$).

Let $\begin{pmatrix} 1 & 1 \\ 0 & 1 \end{pmatrix}$ be a generator of its unipotent subgroup. Let $\phi$ be a cusp form with eigenvalue $\mu = \frac{1-\lambda^2}{4}$. We have then the following Fourier decomposition (see [B]):
$$\phi(x+iy) = \sum_{n \neq 0} a_n y^{\frac{1}{2}} K_{\frac{\lambda}{2}}(2\pi |n|y) e^{2\pi i n x},$$
where $K_{\frac{\lambda}{2}}$ is the $K$-Bessel function.

In order to study the coefficients $a_n$, Rankin and Selberg introduced the series $L(s) = \sum_{n>0} \frac{|a_n|^2}{n^s}$, the Rankin-Selberg $L$-function (we assume that $\phi$ is real valued; hence, $a_n = a_{-n}$). The significance of this Dirichlet series is that it has an integral representation and as a result a spectral interpretation (as well as an analytic continuation!) which we will use.

Let $E(s)$ be the Eisenstein series associated to the cusp at $\infty$. The series $E(s)$ is unitary for $Re(s) = 1/2$ and
$$L(s) = \frac{2\pi^s \Gamma(s)}{\Gamma(s/2)^2 \Gamma(s/2+it)\Gamma(s/2-it)} \langle \phi^2, E(s)\rangle \ ;$$
hence, our method gives an upper bound for $L(s)$. Namely, taking into account the asymptotic behavior of the $\Gamma$-function we obtain, for example, the following:

COROLLARY 1. $\int_T^{T+1} |L(\frac{1}{2}+i\tau)| d\tau \ll T(\ln T)^{\frac{3}{2}}$.

The Lindelöf conjecture for $L(s)$ is stronger: it asserts a bound $|L(\frac{1}{2}+iT)| \ll T^\varepsilon$.

Corollary 1 implies, in turn, a bound on the coefficients $a_n$ themselves via standard methods of analytic number theory (for details, see [G], [P]):

COROLLARY 2. $|a_n| \ll n^{\frac{1}{3}+\varepsilon}$ for any $\varepsilon > 0$.

*Remarks.* 1. The bound $|a_n| \leq cn^{\frac{1}{2}}$ is due to Hecke and follows from the fact that the function $\phi$ is bounded (sometimes this is called the *standard* or *convexity* bound).

The *Peterson-Ramanujan Conjecture* is the assertion that $|a_n| \ll n^{\varepsilon}$ for the congruence subgroups.

The best-known bound for the congruence subgroups is $n^{\frac{5}{28}+\varepsilon}$ due to Bump-Duke-Hoffstein-Iwaniec ([B-I]).

For *nonarithmetic* subgroups, however, there was no improvement over the Hecke bound before [S] appeared. It was even suspected that the Hecke bound might be of true order for nonarithmetic subgroups.

Recently for the general lattice, Sarnak [S] gave the first improvement over the Hecke bound (he treated $SL(2,\mathbb{C})$, while the $SL(2,\mathbb{R})$-case was done in [P]). Sarnak also suggested that the Peterson-Ramanujan Conjecture might be true in this general setting. It was his idea to use the analytic continuation which led us to think about the problem.

2. The main point of Corollary 2 is that it holds without any assumption on the *arithmeticity* of $\Gamma$.

We would like to add that, even theoretically, the triple product method cannot give the Peterson-Ramanujan Conjecture; indeed, even Lindelöf's conjecture for $L(s)$ above implies only that $|a_n| \ll n^{\frac{1}{4}+\varepsilon}$.

The results of this paper where announced in [BR].

*Acknowledgments.* We would like to thank Peter Sarnak for turning our attention to the problem, for fruitful discussions and for initiating our cooperation. We would also like to thank Stephen Semmes for enlightening discussions.

We would like to thank the Binational Science Foundation. Most of the work on this paper was done in a framework of a joint project with P. Sarnak supported by BSF grant No. 94-00312/2.

The second author would like to thank several institutions for providing him with a (temporary) roof while this work was done.

## 1. Analytic continuation of representations

1.1. Let $G$ be a Lie group, $(\pi, G, V)$ its representation and $v$ an analytic vector in $V$. Then we can find a left $G$-invariant domain $U \subset G_{\mathbb{C}}$ containing $G$ such that the function $\xi_v : G \to V$ given by $g \mapsto \pi(g)v$ has an extension to $U$ as a univalued holomorphic function. For the elements $g \in U$ we define the vector $\pi(g)v$ to be the value of the extended function of $\xi_v$ at $g$.

One should be careful with the choice of $U$ since the vector $\pi(g)v$ depends on this choice. However, having fixed $U$, we see that the action of $G$ on $v$ can be unambiguously extended to this somewhat larger set $U \supset G$. We will see that in many situations there is a natural choice of $U$ which works for many vectors $v$.

It is clear that with an appropriate choice of domains of definition the extended operators $\pi(g)$ have the usual properties:

(i) $\pi(gh) = \pi(g)\pi(h)$; $\pi(g^{-1}) = \pi(g)^{-1}$;
(ii) If $\nu : (\pi, V) \to (\tau, L)$ is a morphism of representations, then $\tau(g) \circ \nu = \nu \circ \pi(g)$;
(iii) If $(\omega, V \otimes L)$ is the tensor product of representations $(\pi, V)$ and $(\tau, L)$, then $\omega(g) = \pi(g) \otimes \tau(g)$. If $(\pi^*, V^*)$ is the dual representation, then $\pi^*(g) = \pi(g)^*$.
(iv) If $(\bar\pi, \bar V)$ is the complex conjugate representation, then $\overline{\pi(g)} = \bar\pi(\bar g)$. In particular, given a $G$-invariant positive definite scalar product on V we formally get $\pi(g)^+ = \pi(\bar g)^{-1}$.

1.2. *Geometry of the domain $U$ for* $\mathrm{SL}(2, \mathbb{R})$. (See also Appendix C.) We consider representations of the principal series of the group $G = \mathrm{SL}(2, \mathbb{R})$. Namely, for any $\lambda \in \mathbb{C}$ we consider the representation $(\pi_\lambda, G, \mathfrak{D}_\lambda)$; see 0.1.

In such a realization, the $K$-fixed vector is the function $v(x, y) = (x^2 + y^2)^{\frac{\lambda-1}{2}}$. For convenience, we denote $x^2 + y^2$ by $Q(x, y)$ and will view it as a quadratic form on $\mathbb{C}^2$. Then the action of $G$ on $v$ is given by

$$(1.2) \qquad (\pi(g)v)(x, y) = (g(Q)(x, y))^{(\lambda-1)/2} .$$

Let $U$ be the open subset of $G_\mathbb{C}$ consisting of matrices $g$ such that the quadratic form $g(Q)$ on $\mathbb{R}^2$ has a positive definite real part. Since the function $z \mapsto z^{(\lambda-1)/2}$ is a well-defined holomorphic function in the right half-plane Re $z > 0$, we see that formula (1.2) makes sense for all $g \in U$.

This gives us a holomorphic function on $U$ with values in $\mathfrak{D}_\lambda$. We will see that $U$ is connected, so this function is the holomorphic extension of the function $\xi_v$ to the domain $U$. We will also show that for most $\lambda$ the domain $U$ is the maximal domain of holomorphicity for the function $\xi_v$.

Observe that $U$ is left $G$-invariant and right $K_\mathbb{C}$-invariant, where $K_\mathbb{C} = \mathrm{SO}(2, \mathbb{C}) \simeq \mathbb{C}^* \subset \mathrm{SL}(2, \mathbb{C})$. Let us identify $G_\mathbb{C}/K_\mathbb{C}$ with the variety $\mathcal{Q}$ of unimodular quadratic forms on $\mathbb{C}^2$ via $g \mapsto g(Q)$. By definition, $U$ is the preimage of the open subdomain $\mathcal{Q}_+ \subset \mathcal{Q}$ consisting of all quadratic forms whose real part is positive definite.

For every $\lambda$ we have constructed a holomorphic $G$-equivariant function $v : \mathcal{Q}_+ \to \mathfrak{D}_\lambda$ such that $R \mapsto v_R = R^{(\lambda-1)/2}$, $R \in \mathcal{Q}_+$. The analytic continuation $\pi(g)v$ is given by $\pi(g)v = v_{g(Q)}$.

*Remarks.* 1. Note that all $K$-finite *matrix coefficients* $\langle \pi(g)e_0, e_n \rangle$ have an analytic extension to a much larger domain: $\{\mathrm{diag}(z, z^{-1}) : |\arg(z)| < \frac{\pi}{2}\}$. Observe a curious phenomenon: each matrix coefficient of the function $\pi(g)e_0$ is holomorphic in this larger domain, but the function itself admits analytic continuation to $U$ only.

For groups of higher rank the situation is much more intriguing and we hope to return to it elsewhere.

2. The same proof can be applied to any $K$-finite vector $v \in \mathfrak{D}_\lambda$; it shows that for every such vector the function $\xi_v = \pi(g)v$ has an extension to the same domain $U \subset \mathrm{SL}(2, \mathbb{C})$.

## 2. Triple products. Proof of Theorem 0.2

Recall (see 0.2 and 0.4) that we fix an automorphic function $\phi$ and consider the function $\phi^2 \in L^2(Y) \subset L^2(X)$. Let $\{\phi_i\}$ be the orthonormal eigenbasis of the space $L^2(Y)$, $\Delta \phi_i = \frac{1-\lambda_i^2}{4} \phi_i$. We set $c_i = \langle \phi^2, \phi_i \rangle$, $b_i = |c_i|^2 \exp(\frac{\pi}{2}|\lambda_i|)$. Let $L_i \subset L^2(X)$ be the subspace corresponding to $\phi_i$. We denote by $\mathrm{pr}_i : L^2(X) \to L_i$ the orthogonal projection and by $L^\perp$ the orthogonal complement to the sum of all subspaces $L_i$ in $L^2(X)$.

2.1. *Proof of* (0.4.1). Observe that the Plancherel formula gives us (0.4.1) with an additional term on the right-hand side. The term is equal to $\|\pi(g)(\psi)\|^2$, where $\psi$ is the orthogonal projection of the function $\phi^2$ onto $L^\perp$. Since $L^\perp$ does not have $K$-invariant vectors, $\psi = 0$.

2.2. *Estimates of* $\|\phi_\varepsilon\|$. Choose a family of elements $g_\varepsilon$ tending to the boundary of $U$. Consider the corresponding vectors $v_\varepsilon = \pi(g_\varepsilon)v \in \mathfrak{D}_\lambda$, $v_{i,\varepsilon} = \pi(g_\varepsilon)v_i \in \mathfrak{D}_{\lambda_i}$ and the corresponding functions $\phi_\varepsilon$, $\phi_{i,\varepsilon}$ on $X$. Observe that all our formulas are given not in terms of the element $g_\varepsilon$ (see 0.4) but in terms of the corresponding quadratic form $Q_\varepsilon = g_\varepsilon(Q) \in \mathcal{Q}_+$ (see 1.2). So it is easier for us to describe the forms $Q_\varepsilon$ without specifying elements $g_\varepsilon$.

In our method, the quadratic forms $Q_\varepsilon$ lying within the same $G$-orbit lead to the same estimates; in particular, we can take the diagonal elements $g_\varepsilon$ described in 0.4. Computationally, however, it is easier to work with another system of quadratic forms, namely with the forms $Q_\varepsilon(x, y) = a(x - i\varepsilon y)(\varepsilon x + iy)$, where $i = \sqrt{-1}$ and $a > 0$ is a (bounded as $\varepsilon \to 0$) normalization constant which makes $\det Q_\varepsilon = 1$.

We will see in Appendix C that, modulo the $G$-action, the forms $R \in \mathcal{Q}_+$ depend only on one parameter, so the specific choice of the family $Q_\varepsilon$ is inconsequential.

We can rewrite formula (0.4.2) as

(2.1) $$||\phi_\varepsilon^2||^2 = \sum |c_i|^2 ||\phi_{i,\varepsilon}||^2 \,.$$

PROPOSITION. *Let $(\pi, G, L)$ be an irreducible unitary representation of* $\mathrm{SL}(2,\mathbb{R})$ *and $v \in L$ a unit $K$-fixed vector. Consider $g_\varepsilon$ and $v_\varepsilon = \pi(g_\varepsilon)v$ as above. Then*

(1) $||v_\varepsilon||^2 \leq C|\ln(\varepsilon)|$ *as $\varepsilon \to 0$.*
(2) *There exists $c > 0$ such that if $\pi \simeq \pi_\lambda$ is a representation of the principal series, then $||v_\varepsilon||^2 > c\exp((\frac{\pi}{2} - 6\varepsilon)|\lambda|)$ for any $\lambda$ and $\varepsilon < 0.1$.*
(3) *Fix an isometric $G$-equivariant embedding $\nu : L \to L^2(X)$ and set $\phi_\varepsilon = \nu(v_\varepsilon) \in C^\infty(X)$. Then $\sup_{x \in X} |\phi_\varepsilon(x)| \leq C|\ln \varepsilon|$ as $\varepsilon \to 0$.*

2.3. *Proof of Theorem* 0.2. From Proposition 2.2 it follows immediately that we have:

$$\sup_{x \in X} |\phi_\varepsilon(x)| \leq C|\ln \varepsilon| \text{ and } ||\phi_\varepsilon||^2 = ||v_\varepsilon||^2 \leq C|\ln(\varepsilon)|.$$

Therefore, $||\phi_\varepsilon^2||^2 \leq ||\phi_\varepsilon||^2 \cdot \sup |\phi_\varepsilon|^2 \leq C|\ln(\varepsilon)|^3$. Hence, formula (2.1) implies

$$C|\ln(\varepsilon)|^3 \geq ||\phi_\varepsilon^2||^2 = \sum_i |c_i|^2 ||\phi_{i,\varepsilon}||^2 \geq \sum_i |c_i|^2 e^{(\frac{\pi}{2} - 6\varepsilon)|\lambda_i|} = \sum_i b_i e^{-6\varepsilon|\lambda_i|}.$$

Set $\varepsilon = 1/T$ and collect the terms with $|\lambda_i| \leq T$, and the desired bound results. □

## 3. Invariant norms and estimates of automorphic functions

In this section we prove the upper bound (3) from Proposition 2.

3.1. Let $(\pi, G, L)$ be a unitary representation and $\nu : L \to L^2(X)$ a continuous $G$-equivariant morphism. Then $\nu$ maps the subspace of smooth vectors $V = L^\infty \subset L$ into $C^\infty(X)$. Given a vector $v \in V$, we would like to describe an effective method for obtaining a pointwise bound for the function $\phi = \nu(v)$. In other words, consider the supremum norm $N_{\sup}$ on $V$ defined in 0.5. We would like to find bounds for $N_{\sup}$ in terms of $\pi$.

Observe that the $L^2$-norm of $\phi$ is bounded by $||\nu|| \cdot ||v||$, where $||\nu||$ is the operator norm. So let us assume that $||\nu|| \leq 1$.

First, we will describe some weak bounds of $N_{\sup}$ in terms of Sobolev norms on $V$; these bounds easily follow from the Sobolev restriction lemma. Then we will improve these bounds using the $G$-invariance of $N_{\sup}$.

For convenience we recall the notion of Sobolev norms.

3.2. *Sobolev norms.* Let $(\pi, V)$ be a smooth representation of a Lie group $G$ and $||\cdot||$ be a $G$-invariant Hermitian norm on $V$. For every nonnegative integer $k$ define the *Sobolev norm* $S_k$ on $V$ as follows. Fix a basis $X_1, \ldots, X_n$ of the Lie algebra $\mathfrak{g} = \text{Lie}(G)$ and define the norm $S_k$ by $S_k(v)^2 = \sum ||X_\alpha v||^2$, where the sum runs over all monomials $X_\alpha = X_{i_1} X_{i_2} \cdots X_{i_l}$ of degree $\leq k$.

*Remarks.* 1. If we start with an arbitrary norm $||\cdot||$ on $V$, we get another system of norms, also called Sobolev norms. If the norm $||\cdot||$ is Hermitian, then all Sobolev norms are also Hermitian ones.

Our definition depends on the choice of basis $X_i$ but different choices lead to equivalent norms.

2. Since the norm $||\cdot||$ is $G$-invariant, the representation $(\pi, V)$ is continuous with respect to the norm $S_k$ for any $k$, with continuity constants *independent* of the representation $\pi$. Namely, for every $g \in G$ we have $S_k(\pi(g)v) \leq ||g||_{\text{ad}}^k S_k(v)$, where $||\cdot||_{\text{ad}}$ is the norm in the adjoint representation of $G$.

3. One can actually define Sobolev norm $S_s$ for every $s \in \mathbb{R}$ as follows. The operator $\Delta = -\sum X_i^2 : V \to V$ is an essentially self-adjoint operator on $V$. We can define the Sobolev norm $S_s$ on $V$ to be $S_s(v) = ||(\Delta + 1)^{s/2} v||$.

*Example.* Let $(\pi, V = \mathfrak{D}_\lambda)$ be the unitary representation of the principal series of $G = \text{SL}(2, \mathbb{R})$ and $||\cdot||$ the standard invariant Hermitian norm; $V$ can be identified with $C^\infty_{\text{even}}(S^1)$ and $e_k = e_k(\theta) = e^{2ik\theta}$, $k \in \mathbb{Z}$, is a basis consisting of $K$-finite vectors. For a smooth vector $v$ we define its Fourier coefficients as $a_k = \langle v, e_k \rangle$.

It is easy to check that in this realization the Sobolev norm $S_s$ is the norm induced by the quadratic form $Q_s(v) = \sum_n |a_n|^2 (1 + \mu + 2n^2)^s$ (here we started with any basis of $\mathfrak{g}$ orthonormal with respect to the standard scalar product).

3.3. *Sobolev estimate.* Let $(\pi, G, L)$ be a unitary representation of $G = \text{SL}(2, \mathbb{R})$ and $V \subset L$ the subspace of smooth vectors. Suppose that $X = \Gamma \backslash G$ is compact. Then any morphism of $G$-modules $\nu : V \to C^\infty(X)$ defines the supremum norm $N_{\sup}$ on $V$.

LEMMA 3.1. *Suppose that $||\nu|| \leq 1$ with respect to the $L^2$-norm. Then $N_{\sup} \leq C S_2$, where the constant $C$ only depends on the geometry of $X$.*

The proof of the lemma easily follows from the Sobolev restriction lemma. We will present it in Appendix B together with a similar result for noncocompact lattices.

*Remark.* In [BR] we showed that the same bound $N_{\sup} \ll S_s$ holds for any $s > 1/2$, which is less trivial since it goes beyond the restriction theorem. For our present purposes, however, the elementary result of the lemma is enough.

3.4. *Invariant (semi-)norms.* The bound which we proved in Lemma 3.1 is rather weak. For example, it gives a bound on $N_{\sup}(v_\varepsilon)$ which is a power of $\varepsilon^{-1}$ (even if we use optimal constant $s = 1/2$; see Remark 3.3). We are able to significantly improve this bound using the fact that the norm $N_{\sup}$ is $G$-invariant.

Let us state some elementary general result about invariant (semi-)norms. Let $G$ be an arbitrary group acting on some linear space $V$.

*Claim.* For any seminorm $N$ on $V$ there exists a unique seminorm $N^G$ on $V$ satisfying the following conditions:

(1) $N^G$ is $G$-invariant;

(2) $N^G \leq N$;

(3) $N^G$ is the maximal seminorm satisfying conditions (1) and (2).

We will prove this claim in Appendix A.

The passage from $N$ to $N^G$ has the following obvious properties:

(1) If $N_1 \leq CN_2$, then $N_1^G \leq CN_2^G$;

(2) If $N$ is $G$-invariant, then $N = N^G$.

We apply this general construction to our situation, when the space $V$ is the smooth part of some unitary representation $(\pi, G, L)$ of $G = \mathrm{SL}(2, \mathbb{R})$. Consider the Sobolev norm $S = S_2$ on $V$ and construct the corresponding invariant seminorm $S^G$. If $\nu : L \to L^2(X)$ is a morphism of representations, then $\nu(V) \subset C^\infty(X)$ and we can define the norm $N_{\sup}$ on $V$ as in 0.5. This norm is $G$-invariant and $N_{\sup} \leq CS$. Hence, $N_{\sup} \leq CS^G$; in particular, $N_{\sup}(v_\varepsilon) \leq CS^G(v_\varepsilon)$.

The norm $S^G$, however, is defined in terms of the representation $\pi$ only. It does not depend on the embedding $\nu$. In particular, we can estimate the norm $S^G(v_\varepsilon)$ by computations in $\mathfrak{D}_\lambda$. The main result in this direction is the following proposition which implies inequality (3) in Proposition 2.

PROPOSITION. *Let $(\pi, G, L)$ be a unitary irreducible representation and $v \in L$ a unit $K$-fixed vector. For $k \geq 0$, consider the Sobolev norm $S = S_k$ on the space $V$ of smooth vectors in $L$ and denote by $S^G$ its invariant part.*

*Then there exists a constant $C > 0$ such that $S^G(v_\varepsilon) \leq C|\ln(\varepsilon)|$ as $\varepsilon \to 0$.*

We will prove this proposition in Section 5.

### 4. Noncocompact $\Gamma$

4.1. *Cuspidal representations.* In order to prove the crucial bound, $|\phi_\varepsilon| \leq C|\ln \varepsilon|$, we have used the norm $N_{\sup}$ induced by the supremum norm on $X$ via the embedding $\nu$ and the fact that an appropriate Sobolev norm majorizes it. From this, the proof of the bound and Theorem 0.2 immediately follow. We will explain now how to find such a Sobolev norm in the case of a noncocompact lattice $\Gamma$.

If $X$ is noncompact, it is not clear why a supremum norm exists on the space of smooth vectors of $\pi$. Actually, there is no such norm for a general automorphic representation since a general automorphic function does not need to decay at infinity. However, if $\pi$ is cuspidal, then its smooth vectors decay at infinity and the supremum norm is well defined. A simple proposition below (proven in Appendix B) shows that there is an appropriate Sobolev norm which majors $N_{\sup}$ in this case as in the cocompact case. This suffices to prove the bound $|\phi_\varepsilon| \leq C|\ln \varepsilon|$, hence, the analog of Theorem 0.2.

PROPOSITION. *Let $(\pi, G, L)$ be a unitary representation of the group $G = \mathrm{SL}(2,\mathbb{R})$ and $\nu : L \to L^2(X)$ a bounded morphism of representations whose image lies in the cuspidal part of $L^2(X)$. Consider the space $V = L^\infty$ of smooth vectors in $L$ and introduce the norm $N_{\sup}$ on $V$ as in 0.5. Then there exists a constant $C$ such that $N_{\sup} \leq C S_3$, where $S_3$ is the third Sobolev norm on $V$.*

4.2. We state now the version of Theorem 0.2 for a noncocompact lattice $\Gamma$ (for notations see [B]). Denote by $\{\alpha_j\}_{j=1,\ldots,k}$ the set of cusps and by $E_j(s)$ the corresponding Eisenstein series; let $\{\phi_i\}$ be the basis for the discrete spectrum (cusp forms and residual eigenfunctions). Let $\phi$ be a cusp form and denote, as before, $b_i = |\langle \phi^2, \phi_i\rangle|^2 \exp(\frac{\pi}{2}|\lambda_i|)$ and $b_j(t) = |\langle \phi^2, E_j(\frac{1}{2}+it)\rangle|^2 \exp(\frac{\pi}{2}|t|)$.

THEOREM. *There exists a constant $C$ such that*
$$\sum_{|\lambda_i|\leq T} b_i \ + \ \sum_j \int_{|t|\leq T} b_j(t)dt \leq C(\ln T)^3 \text{ as } T \to \infty.$$

### 5. Some computations in the model $\mathfrak{D}_\lambda$

This section is devoted to the proof of Propositions 2 and 3.4. Our proof is based on explicit computations in the model $\mathfrak{D}_\lambda$ of the representation $\pi$.

Since $\pi$ is a unitary representation with a $K$-fixed vector, it is either a representation of the principal series, or a representation of the complementary series (or the trivial representation). In 5.1 and 5.2 we consider representations of the principal series. In 5.5 we treat the complementary series.

5.1. *Proof of statements* (1) *and* (2) *of Proposition* 2. What we claim in (1) and (2) is independent of the realization of $\pi_\lambda$. We chose the realization of $\pi_\lambda$ in $\mathfrak{D}_\lambda$. By definition, the element $g_\varepsilon$ is chosen so that $v_\varepsilon = \pi(g_\varepsilon)v$ is given by the function $Q_\varepsilon^{\frac{\lambda-1}{2}}$, where $Q_\varepsilon(x,y) = a(x - i\varepsilon y)(\varepsilon x + iy)$.

For computations we will use two models of the representation $\mathfrak{D}_\lambda$:

*Circle model.* Realization of $\mathfrak{D}_\lambda$ as the space of smooth functions on $S^1$, described in 0.1.

*Line model.* In this model, to every vector $v \in \mathfrak{D}_\lambda$ we assign the function $u$ on the line given by $u(x) = v(x, 1)$.

The line model is convenient to describe the action of the Borel subgroup.

LEMMA. (1) $\pi(\begin{pmatrix} 1 & b \\ 0 & 1 \end{pmatrix})u(x) = u(x - b)$.

(2) $\pi(\begin{pmatrix} a & 0 \\ 0 & a^{-1} \end{pmatrix})u(x) = |a|^{\lambda-1}u(a^{-2}x)$.

(3) *For $\lambda = it$ the scalar product in $\mathfrak{D}_\lambda$ is given, up to a factor, by the standard $L^2$-product in the space of functions on the line, namely, $||v||^2 = \frac{1}{\pi}\int |u|^2 dx$.*

Denote by $q_\varepsilon$ the restriction of the quadratic form $Q_\varepsilon$ on the line $\{(x,1)\}$; i.e., $q_\varepsilon(x) = a(x - i\varepsilon)(\varepsilon x + i) = a(\varepsilon(x^2+1) + ix(1-\varepsilon^2))$. Thus, the vector $v_\varepsilon \in \mathfrak{D}_\lambda$ corresponds to the function $u_\varepsilon = q_\varepsilon^{(\lambda-1)/2}$, and we have to estimate the integral $||v_\varepsilon||^2 = \int |u_\varepsilon|^2 dx$.

Let $m(X) = |q(x)|$ and $a(x) = \arg(q(x))$ be the modulus and the argument of the function $q$. Then for $\lambda = it$ we have $|u_\varepsilon(x)|^2 = m(x)^{-1}\exp(2ta(x))$.

*Proof of* (1) *in Proposition* 2. Since $t$ is fixed, the function $\exp(2ta(x))$ is uniformly bounded, while the function $m(x)^{-1}$ is bounded by $\varepsilon^{-1}$ for $|x| \leq \varepsilon$, by $|1/x|$ for $\varepsilon \leq |x| \leq \varepsilon^{-1}$ and by $\varepsilon^{-1}x^{-2}$ for $|x| > \varepsilon^{-1}$, which implies that $\int |u_\varepsilon(x)|^2 dx \leq C|\ln \varepsilon|$. □

*Proof of* (2) *in Proposition* 2. We can assume that $t > 0$. Clearly, on the segment $[1,2]$ we have, uniform in $\varepsilon < 0.1$, bounds $|m(x)| < 3$ and $a(x) > \pi/4 - 3\varepsilon$. This implies that $||v_\varepsilon||^2 \geq c\exp(\pi/2 - 6\varepsilon)$. □

*Remark.* There is another way to compute the norm $||v_\varepsilon||$, based on the theory of spherical functions. Namely, for every $\lambda \in \mathbb{C}$ we consider the spherical function $S_\lambda$ on $G$ equal to the matrix coefficient of the $K$-fixed vector $v \in \mathfrak{D}_\lambda$, $S_\lambda(g) = \langle \pi(g)v, v \rangle$. This function is well-known: it is determined by its restriction to the diagonal subgroup and on this subgroup it is essentially given by the Legendre function. In particular, this function has an analytic continuation to some domain which contains all diagonal matrices $\text{diag}(a^{-1}, a)$ with $|\arg(a)| < \pi/2$.

We can compute the norm $||\pi(g)v||$ using spherical functions as follows. For $g \in U$ we write $||\pi(g)v||^2 = \langle \pi(g)v, \pi(g)v \rangle = \langle \pi(g'g)v, v \rangle = S_\lambda(g'g)$, where $g' = \bar{g}^{-1}$ (see 1.1). In particular, if $g = \text{diag}(a^{-1}, a)$, where $a \in \mathbb{C}$ such that $||a|| = 1$ and $|\arg(a)| < \pi/4$, then we have $||\pi(g)v||^2 = S_\lambda(g^2)$.

5.2. *Proof of Proposition* 3.4. We work with a fixed $\lambda$ as $\varepsilon \to 0$. Denote the norm $S_k^G$ on the space $\mathfrak{D}_\lambda$ by $N$. We want to estimate $N(v_\varepsilon)$.

*Step* 1. The vector $v_\varepsilon$ is realized as the function $Q_\varepsilon^{\lambda-1}$. Consider this function in a circle model. We can choose a partition of unit $\alpha_i$ on the circle and replace the function $v_\varepsilon$ with a function $\alpha v_\varepsilon$, where $\alpha$ is a smooth function with small support on the circle.

If $\alpha$ is supported far from the $x$- and the $y$-axes, then the family of functions $\alpha v_\varepsilon$ is uniformly bounded with respect to the norm $S_k$, hence, with respect to the norm $N$. The case of a function $\alpha$ supported near the $x$-axis can be reduced to the case of the $y$-axis by the change of coordinates $(x \mapsto y, y \mapsto -x)$.

Thus, it suffices to estimate $N(\alpha v_\varepsilon)$, where $\alpha$ is a smooth function supported near the $y$-axis.

*Step* 2. Let us pass to the line model of the representation $\mathfrak{D}_\lambda$. Here one should be a little careful since the standard Sobolev norm $S^k$ on the space $\mathfrak{F}$ of functions on the line does not agree with the Sobolev norm $S^k$ on the space $\mathfrak{D}_\lambda$. However, on the subspace $\mathfrak{F}'$ of functions supported on the segment $[-2, 2]$ these two norms are comparable, and so on this subspace we will pass from one of these norms to another without changing notations.

In the line model our vector $\alpha v_\varepsilon$ is represented by the function $u_\varepsilon$ given by $u_\varepsilon(x) = \alpha a^\kappa (x - i\varepsilon)^\kappa (\varepsilon x + i)^\kappa$, where $\kappa = (\lambda - 1)/2$. We see that as $\varepsilon \to 0$ the structure of the function $u_\varepsilon$ is mainly determined by the factor $(x - i\varepsilon)^\kappa$ which is roughly homogeneous in $x$. We estimate the norm $N(u_\varepsilon)$ using the fact that the norm $N$ itself is homogeneous with respect to dilations. We will do this using the, standard in harmonic analysis, method of dyadic decomposition. Let us describe this method informally for $\lambda = 0$.

In this case, the function $u = u_\varepsilon$ on $[0, 1]$ is, more or less, equal to $(x - i\varepsilon)^{-1/2}$. In other words, $u_\varepsilon$ is just a branch of the function $x^{-1/2}$ slightly smoothed at the origin.

The only *a priori* estimate of the norm $N$ we know is $N \leq S_k$. However, one can easily see that the value $S_k(u)$ is too big. What we can do is to break the segment $I = [0, 1]$ into smaller segments $I_1 = [1/2, 1]$, $I_2 = [1/4, 1/2], \ldots, I_l$ (plus some small segment at the origin) and to break our function $u$ into the sum of functions $u_i$ approximately supported on these segments.

Now let us estimate, separately, the norms $N(u_i)$. The operator $\pi(g)$ with a suitable diagonal matrix $g$ moves $u_i$ into the function $u'_i$ with support on $[1, 2]$. This transformation does not affect the norm $N$, since $N$ is invariant, but it tremendously decreases the Sobolev norm $S_k$. This yields a much better estimate: $N(u_i) = N(u_i') \leq S_k(u_i')$.

To get a better bound, we move the function $u_i$ as far to the right as possible. On the other hand, we cannot move it beyond the point 2 since there we lose control of the Sobolev norm $S_k$; this explains, in particular, why we have to break the function $u$ into pieces: each piece must be scaled differently.

Let us formulate a general statement about functions on the line that sums up the results one can prove using this method.

5.3. *Dyadic decomposition.* Let $\mathfrak{F}$ be the space of smooth functions with compact support on the line. For every $t > 0$ consider the dilation operator $h_t : \mathfrak{F} \to \mathfrak{F}$, where $h_t(f)(x) = f(t^{-1}x)$.

Suppose on $\mathfrak{F}$ we have a homogeneous norm $N$ of degree $r$; i.e., $N(h_t f) = t^{-r} N(f)$. Assume also that for functions supported on the segment $[-2, 2]$ we have the estimate $N(f) \leq S_k(f)$, where $S_k$ is the $k^{\text{th}}$ Sobolev norm.

To estimate the values $N(u_\varepsilon)$ for some family of functions $u_\varepsilon \in \mathfrak{F}$ as $\varepsilon \to 0$, we assume that the family $u_\varepsilon$ is "roughly homogeneous." This means that $u_\varepsilon = \tau_\varepsilon f_\varepsilon \in \mathfrak{F}$, where $f_\varepsilon$ is a family of smooth functions on the line such that $f_{t\varepsilon} = t^\kappa h_t(f_\varepsilon)$; i.e., $f_{t\varepsilon}(tx) = t^\kappa f_\varepsilon(x)$ (we say that this family *is homogeneous of degree $\kappa$*) and $\tau_\varepsilon \in \mathfrak{F}$ is a family of truncation multipliers.

PROPOSITION. *Let $N$ be a norm homogeneous of degree $r$ on the space $\mathfrak{F} = C_c^\infty(\mathbb{R})$. Let $u_\varepsilon \in \mathfrak{F}$ be a family of functions described above. Assume that*:

(1) *There exists a constant $S = S_f$ which bounds the Sobolev norm $S_k$ on the segments $[-2, -1]$ and $[1, 2]$ for all functions $f_\varepsilon$ with $0 < \varepsilon < 1$ and also bounds the Sobolev norm $S_k$ of the function $f_1$ on the segment $[-2, 2]$;*
(2) *The truncation family $\tau_\varepsilon$ is uniformly bounded in $C_c^k[-1, 1]$; i.e., all these functions are supported on the segment $[-1, 1]$ and for all $\varepsilon \leq 1$ all their derivatives up to order $k$ are bounded by some constant $C_{\text{tr}}$.*

*Then $N(u_\varepsilon) \leq CC_{\text{tr}} S_f (\varepsilon^{\operatorname{Re}\kappa - r} + \int_\varepsilon^1 t^{\operatorname{Re}\kappa - r} \cdot dt/t)$.*

In other words, $N(u_\varepsilon) \ll 1$ if $\operatorname{Re}\kappa > r$, $N(u_\varepsilon) \ll \varepsilon^{\operatorname{Re}\kappa - r}$ if $\operatorname{Re}\kappa < r$ and $N(u_\varepsilon) \ll |\ln \varepsilon|$ if $\operatorname{Re}\kappa = r$.

We can apply this proposition to our situation. Namely, consider the family of functions $f_\varepsilon(x) = (x + i\varepsilon)^\kappa$, where $\kappa = (\lambda - 1)/2$. Identify the space $\mathfrak{F} = C_c^\infty(\mathbb{R})$ with a subspace in $\mathfrak{D}_\lambda$ using the line model of $\mathfrak{D}_\lambda$ (see 5.1). Then the formulas for the action of the diagonal group on $\mathfrak{F}$ from Lemma 5.1 show that the $G$-invariant norm $N$ on $\mathfrak{D}_\lambda$ considered as a norm on $\mathfrak{F}$ is homogeneous of degree $r = -1/2$.

It is easy to check that the family of functions $u_\varepsilon$ (for some $\tau_\varepsilon$) satisfies the conditions of Proposition 5.3 with $\kappa = (\lambda - 1)/2$. Thus, $\operatorname{Re}\kappa = -1/2$ and Proposition 5.3 shows that $N(u_\varepsilon) \leq C|\ln\varepsilon|$, which proves Proposition 3.4. □

5.4. *Proof of Proposition* 5.3. Let us formalize the proof outlined in 5.2. Technically, it is a little easier to break $u$ into an integral, rather than a sum, of components. That is what we are going to do.

Fix a smooth function $\gamma \in \mathfrak{F}$ equal to 1 on $[-1, 1]$ and supported on $[-2, 2]$. Then $\gamma u_\varepsilon = u_\varepsilon$.

Consider two families of functions in $\mathfrak{F}$: $\gamma_t = h_t(\gamma)$ and $\delta_t = t\frac{d}{dt}\gamma_t$.

They have the following properties:

(i) $\gamma_1$ is supported on $[-2, 2]$; $\delta_1$ is supported on $[-2, -1] \cup [1, 2]$;
(ii) $h_t(\gamma_1) = \gamma_t$, $h_t(\delta_1) = \delta_t$;
(iii) $\gamma_1 - \gamma_a = \int_a^1 \delta_t \cdot dt/t$.

Fix an $\varepsilon \leq 1$ and consider two families of functions in $\mathfrak{F}$: $g_t = \gamma_t u_\varepsilon$ and $c_t = \delta_t u_\varepsilon$. Then we can express $u_\varepsilon = g_1$ as $g_1 = g_\varepsilon + \int_\varepsilon^1 c_t \cdot dt/t$. Hence, the following claim (with $C_{\mathrm{tr}}$, $S_f$ described in the statement of Proposition 5.3) implies the bound on the norm $N(u_\varepsilon)$:

*Claim.* There exists $C$ which only depends on $k$ and such that

(1) $N(c_t) \leq CC_{\mathrm{tr}}S_f t^{\operatorname{Re}\kappa - r}$ for all $t$ such that $\varepsilon \leq t \leq 1$;
(2) $N(g_\varepsilon) \leq CC_{\mathrm{tr}}S_f \varepsilon^{\operatorname{Re}\kappa - r}$.

*Proof.* (1) Set $a = t^{-1}$. Then we have $N(c_t) = t^{-r}N(h_a(c_t))$ and $h_a(c_t) = h_a(\delta_t)h_a(\tau_\varepsilon)h_a(f_\varepsilon) = \delta_1 h_a(\tau_\varepsilon)f_{a\varepsilon} \cdot t^\kappa$. Observe that the function $h_a(c_t)$ is supported on $[-2, -1] \cup [1, 2]$, where all the derivatives of $\delta_1 h_a(\tau_\varepsilon)$ up to the order $k$ are uniformly bounded by $CC_{\mathrm{tr}}$, while the $S_k$-norm of all functions $f_{a\varepsilon}$ is bounded by $S_f$. This implies that $N(h_a(c_t)) \leq S_k(h_a(c_t)) \leq CC_{\mathrm{tr}}S_f t^{\operatorname{Re}\kappa}$. Hence, $N(c_t) \leq CC_{\mathrm{tr}}S_f t^{\operatorname{Re}\kappa - r}$.

Claim (2) is similarly proved with the help of the dilation $h_\varepsilon^{-1}$. □

5.5. *The complementary series.* We describe modifications in the proofs in 5.1 and 5.2 needed to treat the complementary series.

Let $(\pi, G, V)$ be a representation of the complementary series. We will realize it as a representation $(\pi_\lambda, \mathfrak{D}_\lambda)$ for some $\lambda$ such that $-1 < \lambda < 0$ (see formulas in [G5]).

We can use the line model of the space $\mathfrak{D}_\lambda$ as in 5.1 but with the scalar product given by $||f||^2 = \int |x - x'|^{-\lambda - 1}f(x)\overline{f(x')}dxdx'$. As in 5.1, we have to estimate $||u_\varepsilon||^2$, where $u_\varepsilon = (q_\varepsilon)^{(\lambda-1)/2}$. In this case the main contribution comes from a neighborhood of 0 and direct computations show that $||u_\varepsilon||^2 \leq C|\ln\varepsilon|$.

Proof in 5.2 is even easier. Indeed, the invariant norm $N$ on $\mathfrak{D}_\lambda$ becomes a homogeneous of degree $(\lambda - 1)/2$ norm $N$ on the space $\mathfrak{F}$ of functions on the line. We have to estimate the value of this norm on a function $u_\varepsilon$, roughly homogeneous of degree $\kappa = (\lambda - 1)/2$. Proposition 5.3 implies that $N(u_\varepsilon) \leq C|\ln \varepsilon|$. □

## Appendix A. Invariant norms

A.1. Fix a complex vector space $V$. A *seminorm* on $V$ is a function $N : V \to \mathbb{R}^+$ such that $N(v_1 + v_2) \leq N(v_1) + N(v_2)$ and $N(av) = |a|N(v)$.

The set $\mathcal{N}$ of all seminorms on $V$ is a partially ordered set with respect to the relation $\leq$, where $N \leq N'$ whenever $N(v) \leq N'(v)$ for all $v \in V$.

LEMMA. *The partially ordered set $(\mathcal{N}, \leq)$ is inductive, i.e., any nonempty family of seminorms $N_u$ has the exact lower bound $N = \inf N_u$.*

*Proof.* Define $N$ by setting $N(v) = \inf(\sum N_{u_i}(v_i))$, where the sum runs over all finite collections $v_1, \ldots, v_k \in V$ and $N_{u_1}, \ldots, N_{u_k} \in \mathcal{N}$ which satisfy $\sum v_i = v$.

It is easy to check that $N$ is a seminorm bounded by each seminorm $N_u$. If $M$ is a seminorm such that $M \leq N_u$ for all $N_u$, then, clearly, $M \leq N$. The uniqueness of $N$ is obvious. □

Geometrically it is clear that the seminorm $N$ is defined by the unit ball which is the convex hull of the unit balls of the seminorms $N_u$.

A.2. *Construction of invariant seminorms.* Suppose an arbitrary group $G$ acts on a complex vector space $V$. Then $G$ acts on the set of seminorms on $V$ by $g(N)(v) = N(g^{-1}v)$. For every seminorm $N$ on $V$ we define a new seminorm: $N^G = \inf\limits_{g \in G} g(N)$. From the definition one can immediately deduce that:

(i) $N^G$ is an invariant seminorm.
(ii) If $M$ is any invariant seminorm bounded by $N$, then $M \leq N^G$.

*Remark.* Let us apply this construction to a representation $(\pi, G, V)$ of the unitary principal series of the group $G = \mathrm{SL}(2, \mathbb{R})$.

For any $s \in \mathbb{R}$, consider the corresponding Sobolev norm $S_s$ and construct the invariant seminorm $S_s^G$. One can easily show that for $s < 0$ the seminorm $S_s^G$ vanishes. For $s \geq 0$ the norm $S_s$ is bounded below by an invariant unitary norm $||\cdot||$. This implies that $||\cdot|| \leq S_s^G$; hence, $S_s^G$ is a norm.

One can show that for all $s > 1/2$ the norms $S_s^G$ are equivalent to the same norm, $B$. This norm is distinguished by the condition that it is the maximal $G$-invariant norm on $V$ (i.e., any $G$-invariant norm $N$ is bounded by $CB$ for some constant $C$). If we realize $V$ as the space of functions on the circle (see 0.1), then $B$ becomes equivalent to the Besov norm $B_{1,1}^{1/2}$ (here $1/2$ stands for derivative of order $1/2$, one 1 stands for the $L^1$-norm and another 1 is the weight index in the Besov norm). This equivalence can be shown by standard methods of harmonic analysis, like in [A], where similar questions are discussed.

In particular, this shows that for any cuspidal representation we have a bound on the supremum norm of the form $N_{\sup} \ll B$.

For $0 \leq s \leq 1/2$ the norms $S_s^G$ are all distinct. It is interesting to investigate the nature of these norms. We think that in the circle model they are close to Besov norms, namely $S_s^G \sim B_{q,q}^s$ with $1/q - s = 1/2$.

The embedding $\nu : V \to C^\infty(X)$ defines a family of $G$-invariant norms $N_p$ on the space $V$ corresponding to $L^p$-norms on $X$. We can use Proposition 4.1 and the Besov norm $B$ described above to give some bounds for these norms in terms of the representation $(\pi, G, V)$ only.

Namely, consider a representation of the principal series. Realize the space $V$ in the circle model. Then we have the required bounds at two points:

(i) $N_2$ coincides with the $L^2$-norm on $S^1$, the latter being the Besov norm $B_{2,2}^0$.

(ii) $N_\infty$ is bounded by the norm $B = B_{1,1}^{1/2}$.

Using interpolation theory for Banach norms we conclude that for intermediate $p$, $2 \leq p \leq \infty$, there exists a bound $N_p \ll B_{q,q}^s$, where $1/q + 1/p = 1$, $1/q - s = 1/2$.

It is interesting to give similar bounds for norms $N_p$ for $1 \leq p < 2$.

### Appendix B. Estimates using Sobolev norms

B.1. *Sobolev inequalities.* We start with the standard Sobolev lemma:

LEMMA. *Let $B$ be the unit ball in $\mathbb{R}^n$. Consider the space $V$ of smooth functions on $B$ with the $L^2$-norm $||\cdot||$, and introduce Sobolev norms $S_k$ on $V$ as in subsection 3.2, i.e., $S_k(f)^2 = \sum ||\partial_\alpha(f)||^2$, where the sum runs over all monomials in partial derivatives $\partial_\alpha = \partial_{i_1}, \ldots, \partial_{i_l}$ of order $\leq k$. Then for any $k > n/2$ there exists a constant $C$ such that $|f(0)| \leq C S_k(f)$ for $f \in V$.*

This lemma holds (though its formulation is more cumbersome) for any Sobolev norm $S_s$ with $s > n/2$.

In this paper we actually need this lemma only for $k \geq n$. In this case the estimate is elementary and can be proven by induction using direct integration.

From the Sobolev lemma we immediately deduce its version for any Lie group $G$. Namely, suppose $G$ is an $n$-dimensional Lie group. We fix some basis $X_i$ of the Lie algebra $\mathfrak{g}$ of the group $G$ and use it to construct a left-invariant metric and a left-invariant measure on $G$. Fix a symmetric compact neighborhood $B$ of the unit $e \in G$; we will call $B$ the *unit ball*.

Consider the space $V$ of smooth functions on $B$, introduce the $L^2$-norm with respect to the left-invariant measure on $G$ and define Sobolev norms on $V$ using derivatives $X_i \in \mathfrak{g}$ corresponding to the right $G$-action.

COROLLARY. *Let $k > n/2$. Then there exists a constant $C$ such that $|f(e)| \leq CS_k(f)$, for any $f \in V$.*

Though this reformulation seems to be quite trivial, it is in fact rather strong, since it gives some estimates which are *uniform* with respect to the left $G$-action.

B.2. *Sobolev inequalities for homogeneous spaces.* Let $X = \Gamma\backslash G$ be a homogeneous space. We consider a measure on $X$ induced by a left $G$-invariant Haar measure on $G$ and introduce the $L^2$-norm and Sobolev norms in the space $V$ of smooth functions on $X$. We would like to describe Sobolev type inequalities which are uniform on $X$.

We will call a function $w$ on the space $X$ a *weight* if it is a positive measurable function and for every $g \in G$ there exists a constant $C(g)$ such that $w(xg) \leq C(g)w(x)$ for all $x \in X$; we also assume that the function $C(g)$ is locally bounded.

Fix a ball $B \subset G$ as above. For every point $x \in X$ we consider the map $p_x : B \to X$ given by $g \mapsto xg$. This map induces the morphism of Hilbert spaces $p_x^* : L^2(X) \to L^2(B)$ and we denote by $w(x)$ the norm of this morphism. (One can show that $w(x)^2$ is the maximal cardinality of the fibers of the map $p_x$.) It is easy to see that $w$ is a weight on $X$. It depends on the choice of the unit ball $B$ but for different balls these functions are comparable.

Now the Sobolev inequality for $G$ (Corollary B.1) immediately implies

PROPOSITION. *For $k > n/2$ there exists a constant $C$ such that for any $f \in V$ and $x \in X$, $|f(x)| \leq Cw(x)S_k(f)$.*

If $X$ is compact, then the weight function $w(x)$ is bounded; thus the proposition implies Lemma 3.3.

If $X$ is not compact, then in order to get a bound on $\sup |f(x)|$ we need an additional information, e.g., that $f$ is cuspidal.

B.3. *Sobolev estimates in the cuspidal case.* Let $G$ be a Lie group, $\Gamma \subset G$ a lattice in $G$ and $X = \Gamma\backslash G$. We fix a basis in $\mathfrak{g} = \mathrm{Lie}(G)$ and use it to construct a left-invariant metric on $G$ and the induced metric on $X$.

Recall the notion of the cuspidal function on $X$. A unipotent subgroup $U \subset G$ is *cuspidal* if $U$ is nontrivial and $\Gamma_U = \Gamma \cap U$ is a cocompact subgroup in $U$. Geometrically this means that the left orbits of $U$ in $G$ become compact when mapped to $X$; these compact sets are called *horocycles*.

A smooth function $f$ on $X$ is called *cuspidal* if the integral of $f$ over any horocycle vanishes. Important here is the $G$-invariance of the space of cuspidal functions on $X$; in particular, this space is $\mathfrak{g}$-invariant.

For every point $x \in X$, define $d(x)$ as the infimum of 1 and the diameters of all horocycles through $x$. Roughly speaking, $d$ measures how close point $x$ is to a cusp.

The following lemma shows that we can improve estimates for a cuspidal function if we know estimates for its derivatives.

LEMMA. *Suppose $f$ is a smooth cuspidal function on $X$ such that $f$ and all its derivatives $X_i f$ are bounded by a weight $w$. Then for some constant $C$ independent of $f$, the function $f$ is bounded by a weight $w' = Cwd$.*

*Proof.* Let $x \in X$. If $d(x) \geq 0.1$, there is nothing to prove. So suppose $d(x) < 0.1$; i.e., there exists an horocycle $H$ passing through $x$ and $\mathrm{diam}(H) \leq 2d(x) < 0.2$. By hypothesis, all the derivatives of $f$ at all points of $H$ are bounded by $Cw(x)$. This means that the gradient of $f$ is bounded by $C'w(x)$ at all points of $H$.

We can assume that the function $f$ is real. Then the condition that its integral over $H$ vanishes implies that $f$ vanishes at least at one of the points of $H$. Combined with the estimate $|\mathrm{grad}(f)| \leq C'w(x)$ and the fact that $\mathrm{diam}(H) \leq 2d(x)$ this implies the desired estimate $|f(x)| \leq 2C'w(x)d(x)$. □

*Proof of Proposition* 4.1. Let us come back to the case $G = \mathrm{SL}(2, \mathbb{R})$, $X = \Gamma \backslash G$.

Let $f$ be a cuspidal function on $X$. It follows from the estimate in B.2 that $f$ and its derivatives are bounded by $Cw(x)S_3(f)$ with the weight $w$ introduced in B.2. Hence, $|f| \leq Cw(x)d(x)S_3(f)$ by Lemma B.3. In order to finish the proof, it suffices to show that the function $wd$ is bounded on $X$. This result easily follows from the theory of Siegel domains. Let us recall this theory.

Fix a cuspidal unipotent subgroup $U$, set $\Gamma_U = \Gamma \cap U$ and consider the homogeneous space $X_U = \Gamma_U \backslash G$. Let $p : X_U \to X$ be the natural projection.

Let $A$ be the Cartan group of $G$. It is canonically isomorphic to $\mathbb{R}^*$; we define the function $h : A \to \mathbb{R}$ by $h(a) = \ln(|a|)$. Using Iwasawa decomposition we can canonically extend $h$ to a left $U$-invariant and right $K$-invariant function on $G$. This latter function defines the function $h : X_U \to \mathbb{R}$.

A Siegel domain $\mathfrak{S}_T$ is an open subdomain of $X_U$ defined as the preimage $h^{-1}(T, \infty)$. Fix one Siegel domain $\mathfrak{S} = \mathfrak{S}_T$ and consider the map $p : \mathfrak{S} \to X$. It is easy to see that for the unit ball $B \subset G$ the domain $\mathfrak{S}' = \mathfrak{S} \cdot B$ is contained in another Siegel domain, $\mathfrak{S}_{T'}$. By the theory of Siegel domains the cardinality of the fibers of the map $p : \mathfrak{S}' \to X$ is finite and bounded; hence, the operator $p^* : L^2(X) \to L^2(\mathfrak{S}')$ is bounded. By the definition of $w$ this implies that $w(p(x)) \leq Cw(x)$ for any $x \in \mathfrak{S}$. Also, obviously, $d(p(x)) \leq d(x)$ for any $x \in \mathfrak{S}$.

By reduction theory the space $X$ can be covered by a compact set and the union of images of a finite collection of Siegel domains. Hence, in order to check that the function $dw$ is bounded on $X$, it suffices to check the bound on each Siegel domain.

Direct computations on the Siegel domain $\mathfrak{S}_T$ show that for $x \in \mathfrak{S}_T \subset X_U$ we have $d(x) \leq C \exp(-h(x))$ and $w(x) \leq C \exp(h(x)/2)$, which shows that the function $dw$ is bounded on $\mathfrak{S}_T$. $\square$

## Appendix C. The geometry of domain $U \subset \mathrm{SL}(2, \mathbb{C})$

In the 2-dimensional case there is a convenient way to describe the variety $\mathcal{Q}$ of unimodular quadratic forms; see 1.2. Namely, fix a skew-symmetric form $\langle \cdot, \cdot \rangle$ on $\mathbb{C}^2$. Then for any pair of noncollinear vectors $a, b \in \mathbb{C}^2$ set $Q_{a,b}(x) = \frac{\langle a, x \rangle \langle b, x \rangle}{\langle a, b \rangle}$. Note that this form only depends on the images of $a$ and $b$ in $\mathbb{P}^1$, i.e., we can define the unimodular quadratic form $Q_{a,b}$ for any pair of distinct points of $\mathbb{P}^1$. It is easy to see that this identification defines a natural $\mathrm{SL}(2, \mathbb{C})$-equivariant isomorphism of algebraic varieties $\mathbb{P}^1 \times \mathbb{P}^1 \setminus \{\text{diagonal}\}$ and $\mathcal{Q}$.

LEMMA. (1) *Identify $\mathbb{C}$ with an open subset of $\mathbb{P}^1$, $z \mapsto (z, 1)$. (In particular, this realizes $\mathfrak{H}$ as an open $G$-invariant subset of $\mathbb{P}^1$.) Then the morphism $(a, b) \mapsto Q_{a,b}$ identifies the domain $\mathfrak{H} \times \bar{\mathfrak{H}}$ with the subdomain $\mathcal{Q}_+ \subset \mathcal{Q}$.*

(2) *Let $I(Q) = \{\mathrm{diag}(z^{-1}, z) : |\arg z| < \frac{\pi}{4}\} \subset G_\mathbb{C}$. Then $U = \mathrm{SL}(2, \mathbb{R}) \cdot I \cdot \mathrm{SO}(2, \mathbb{C})$.*

*Proof.* If $P$ is a diagonal unimodular quadratic form on $\mathbb{C}^2$, then it is easy to see that $P \in \mathcal{Q}_+$ if and only if $P$ is of the form $g(Q)$ for some $g \in I$. This shows that $I \subset U$ and, hence, $G \cdot I \cdot K_\mathbb{C} \subset U$.

In order to prove the opposite inclusion it suffices to show that $\mathcal{Q}_+ = G \cdot I(Q)$. Consider any form $P \in \mathcal{Q}_+$. Since the form $\mathrm{Re}\, P$ is positive definite, we can find a basis of $\mathbb{R}^2$ in which the real and imaginary parts of $P$ are simultaneously diagonalized. This implies that the orbit $G(P)$ contains a diagonal quadratic form $P'$. As we saw, $P'$ is of the form $h(Q)$ for some $h \in I$. This proves heading (2).

The form $Q$ corresponds to the point $(i, -i) \in \mathfrak{H} \times \bar{\mathfrak{H}}$. Hence, for $h \in I$ the point $(hi, -hi) \in \mathfrak{H} \times \bar{\mathfrak{H}}$ corresponds to the quadratic form $h(Q)$ which lies in $\mathcal{Q}_+$; thus, the subset $I(Q)$ lies in both $\mathfrak{H} \times \bar{\mathfrak{H}}$ and $\mathcal{Q}_+$.

Since $\mathcal{Q}_+$ and $\mathfrak{H} \times \bar{\mathfrak{H}}$ are $G$-invariant, to show that they coincide, it suffices to show that each of them is generated by $I(Q)$ under the $G$-action.

For $\mathcal{Q}_+$ this fact was already proven in heading (2). In order to prove it for $\mathfrak{H} \times \bar{\mathfrak{H}}$, observe that the only $G$-invariant of the pair $(a,b) \in \mathfrak{H} \times \bar{\mathfrak{H}}$ is the hyperbolic distance $d(a,\bar{b})$. Since, in $I(Q)$, arbitrary distances can be realized, $I(Q)$ generates $\mathfrak{H} \times \bar{\mathfrak{H}}$ as a $G$-set. □

Tel Aviv University, Ramat Aviv, Israel
*E-mail address*: bernstei@math.tau.ac.il

Ben Gurion University, Beer Sheva and Weizmann Institute of Science, Rehovot, Israel
*E-mail address*: reznik@wisdom.weizmann.ac.il